\documentclass{amsart}

\usepackage{amsmath} %I used the environment gather once.
\usepackage{amssymb}
\usepackage[final]{epsfig}
\usepackage{hyperref}
\usepackage{xy}
\usepackage{xspace}

\xyoption{all}
\CompileMatrices

\newtheorem{PROP}{Proposition}[section]
\newtheorem{THM}[PROP]{Theorem}
\newtheorem{LM}[PROP]{Lemma}
\newtheorem{COR}[PROP]{Corollary}
\theoremstyle{definition}
\newtheorem{DEF}[PROP]{Definition}
\newtheorem{EXA}[PROP]{Example}
\newtheorem{REM}[PROP]{Remark}

\newcommand{\Kat}[1]{\ensuremath{\mathbf{#1}}}
\newcommand{\Top}{\Kat{Top}}
\newcommand{\FG}{\Kat{FG}}

\newcommand{\Seq}{\Kat{Seq}}
\newcommand{\SSeq}{\Kat{SSeq}}
\newcommand{\Topal}{\Top(\alpha)}
\newcommand{\Topom}{\Top(\omega_1)}
\newcommand{\Genal}{\Kat{Gen}(\alpha)}
\newcommand{\Gen}{\Kat{Gen}}
\newcommand{\CHb}[1]{\ensuremath{\mathrm{CH(}#1\mathrm{)}}}
\newcommand{\CH}[1]{\CHb{\Kat{#1}}}
\newcommand{\CHA}{\CHb A}
\newcommand{\SCHb}[1]{\ensuremath{\mathrm{SCH(}#1\mathrm{)}}}
\newcommand{\SCHA}{\SCHb A}
\newcommand{\SSb}[1]{\ensuremath{\mathrm S #1}\xspace}
\newcommand{\SSp}[1]{\SSb{\Kat{#1}}}
\newcommand{\HCKb}[1]{\ensuremath{\mathrm{HCK(}#1\mathrm{)}}}
\newcommand{\HCK}[1]{\HCKb{\Kat{#1}}}
\newcommand{\filsum}[1]{\ensuremath{#1}-sum\xspace}
\newcommand{\filsums}[1]{\ensuremath{#1}-sums\xspace}
\newcommand{\SumA}[3]{\ensuremath{\sum\limits_{#1} \langle #2, #3 \rangle}}
\newcommand{\Suma}[2]{\ensuremath{\sum \langle #1,#2 \rangle}}
\newcommand{\SuM}[2]{\ensuremath{\sum\limits_{#1} #2}}
\newcommand{\Sum}[1]{\ensuremath{\sum #1}}

\newcommand{\Com}{\ensuremath{C(\omega_0)}}
\newcommand{\Calp}{\ensuremath{C(\alpha)}}
\newcommand{\Cbet}{\ensuremath{C(\beta)}}
\newcommand{\Som}{\ensuremath{A_\omega}}

\newcommand{\On}{\ensuremath{\mathrm{ON}}}
\newcommand{\N}{\ensuremath{\mathbb{N}}}
\newcommand{\emps}{\ensuremath{\emptyset}}

\newcommand{\Zobr}[3]{\ensuremath{#1\colon #2\to #3}}
\newcommand{\Vloz}[3]{\ensuremath{#1\colon #2\hookrightarrow #3}}
\newcommand{\card}[1]{\ensuremath{\operatorname{card}#1}}

\newcommand{\inv}[1]{\ensuremath{#1^{-1}}}
\newcommand{\Invobr}[2]{\inv{#1}[#2]}
\newcommand{\Obr}[2]{\ensuremath{#1[#2]}}

\newcommand{\ol}[1]{\ensuremath{\overline{#1}}}
\newcommand{\Uzav}[2]{\ensuremath{\overline{#1}\mathstrut^{#2}}}

\newcommand{\iaoi}{if and only if\xspace}
\newcommand{\ie}{i.e.\xspace}
\newcommand{\Ie}{I.e.\xspace}
\newcommand{\eg}{e.g.\xspace}
\newcommand{\Wlog}{W.l.o.g.\xspace}
\newcommand{\wlogg}{w.l.o.g.\xspace}
\newcommand{\uv}[1]{``#1''}
\newcommand{\ssp}{subspace\xspace}
\newcommand{\ssps}{subspaces\xspace}
\newcommand{\tsp}{topological space\xspace}
\newcommand{\tsps}{topological spaces\xspace}

\newcommand{\filsp}{prime space\xspace}
\newcommand{\filsps}{prime spaces\xspace}
\newcommand{\filssp}{prime subspace\xspace}
\newcommand{\filssps}{prime subspaces\xspace}
\newcommand{\primfac}{prime factor\xspace}

\newcommand{\primfacs}{prime factors\xspace}
\newcommand{\resp}{resp.\xspace}
\newcommand{\vp}{\ensuremath{\varphi}}
\newcommand{\omnul}{\ensuremath{\omega_0}}
\newcommand{\Sierp}{Sierpi\'nski\xspace}
\newcommand{\rcn}{regular cardinal\xspace}
\newcommand{\arcn}{a regular cardinal\xspace}

\newcommand{\Anm}{\ensuremath{P(A_n)}}
\newcommand{\Aom}{\ensuremath{A_\omega}}

\newcommand{\Can}{\ensuremath{\Calp_n}}
\newcommand{\Cam}{\ensuremath{\Calp_m}}
\newcommand{\Camj}{\ensuremath{\Calp_{m+1}}}
\newcommand{\Canm}{\ensuremath{\Calp_n^-}}

\newcommand{\Canjm}{\ensuremath{\Calp_{n+1}^-}}
\newcommand{\Canom}{\ensuremath{\Calp_{n-1}^-}}
\newcommand{\Camjm}{\ensuremath{\Calp_{m+1}^-}}
\newcommand{\Camm}{\ensuremath{\Calp_m^-}}
\newcommand{\Camk}{\ensuremath{\Calp_m^\eta}}

\newcommand{\TSSd}[1]{TSS'_{#1}(\Com)}
\newcommand{\ciTSS}[1]{T'SS'_{#1}}

\begin{document}

\title{On hereditary coreflective subcategories of $\Top$}
\author{Martin Sleziak}
\address{Department of Algebra and Number Theory, FMFI UK, Mlynsk\'a dolina, 842 48 Bratislava, Slovakia}
\email{sleziak@fmph.uniba.sk}

\maketitle

\begin{abstract}
Let $A$ be a topological space which is not finitely generated and $\CHA$
denote the coreflective hull of $A$ in $\Top$. We construct a generator of the
coreflective subcategory $\SCHA$ consisting of all \ssps of spaces from $\CHA$
which is a \filsp and has the same cardinality as $A$. We also show that if
$\Kat A$ and $\Kat B$ are coreflective subcategories of $\Top$ such that the
hereditary coreflective kernel of each of them is the subcategory $\FG$ of all
finitely generated spaces, then the hereditary coreflective kernel of their
join $\CHb{\Kat A\cup \Kat B}$ is again $\FG$.

\noindent Keywords: coreflective subcategory, hereditary coreflective subcategory,
hereditary coreflective hull, hereditary coreflective kernel, prime space

\noindent 2000 MSC classification: Primary 54B30; Secondary 18B30.
%54B30 Categorical methods (GENERAL TOPOLOGY/Basic constructions)
%18B30 Categories of topological spaces and continuous mappings (CATEGORY THEORY/special categories)
\end{abstract}

\section*{Introduction}

Let $X$ be a \tsp which is not finitely generated and $\SCHb X$ be the
hereditary coreflective hull of $X$ in the category $\Top$ of \tsps. The aim of
this paper is to construct a \filsp $Y_X$ with the same cardinality as $X$ such
that $\SCHb X=\CHb{Y_X}$ where $\CHb{Y_X}$ is the coreflective hull of $Y_X$.
Obviously, if $X$ is finitely generated, then $\CHb X=\SCHb X$. If $X$ is not
finitely generated, then, using the \primfacs of $X$ we can easily construct a
\filsp $P_X$ such that $\SCHb X=\SCHb{P_X}$. Thus, it suffices to restrict our
investigation to the case of \filsps.

For the \filsp $\Com$ consisting of a convergent sequence and its limit point
the problem was studied in \cite{FRARAJ}, where a countable generator for the
category $\SCHb\Com$ of subsequential spaces was produced.

Our procedure of constructing a generator $Y_A$ of the category $\SCHA$ (where
$A$ is a \filsp that is not finitely generated) consists of two main steps. In
the first step, using similar methods as in \cite{FRARAJ}, we produce a set of
special \filsps which generates $\SCHA$. Then, in the second step, we construct
the generator $Y_A$ of $\SCHA$ with the required properties.

This construction was inspired by the space $S_\omega$ from \cite{ARHFRA} and
in the case $A=\Com$ it gives a countable generator for the category of
subsequential spaces different from that one presented in \cite{FRARAJ}.

Finally, as an application of some above mentioned results we prove that if
$\Kat A$ and $\Kat B$ are coreflective subcategories of $\Top$ such that the
hereditary coreflective kernel of $\Kat A$ as well as the hereditary
coreflective kernel of $\Kat B$ is the category $\FG$ of finitely generated
spaces, then $\FG$ is also the hereditary coreflective kernel of their join
$\CHb{\Kat A\cup\Kat B}$. As a consequence of this result and some results of
\cite{SLEZ} we obtain that the collection of all those coreflective
subcategories of $\Top$ the hereditary coreflective kernel of which is $\FG$
and the hereditary coreflective hull of which is $\Top$ is closed under the
formation of non-empty finite joins (in the lattice of all coreflective
subcategories of $\Top$) and arbitrary non-empty intersections.

\section{Preliminaries}

We recall some known facts about coreflective subcategories of the category
$\Top$ of \tsps (see \cite{HER}). All subcategories are supposed to be full and
isomorphism-closed. The topological sum is denoted by $\sqcup$.

Let $\Kat A$ be a subcategory of $\Top$. $\Kat A$ is \emph{coreflective} \iaoi
it is closed under the formation  of topological sums and quotient spaces. If
$\Kat A $ is a subcategory of $\Top$ or a class of \tsps, then the
\emph{coreflective hull} of $\Kat A$ is the smallest coreflective subcategory
of $\Top$ which contains $\Kat A$ and we denote it by $\CH A$. $\CH A$ consists
of all quotients of topological sums of spaces that belong to $\Kat A$. If
$\Kat B= \CH A$, then we say that $\Kat A$ \emph{generates} $\Kat B$ and the
members of $\Kat A$ are called \emph{generators} of $\Kat B$. If $\Kat
B=\CHb{\{X\}}$, then $\Kat B$ is called \emph{simple generated} and $X$ is said
to be a \emph{generator} of $\Kat B$. We use the notation $\Kat B = \CHb X$ in
this case.

Let $\Kat A$ be a subcategory of $\Top$ and let  $\SSp A$ denote the
subcategory of $\Top$ consisting of all subspaces of spaces from $\Kat A$. Then
the following result is known (see \cite[Remark 2.4.4(5)]{KAN} or
\cite[Proposition 3.1]{CIN}).

\begin{PROP}
If \Kat A is a coreflective subcategory of $\Top$, then $\SSp A$
is also a coreflective subcategory of $\Top$. ($\SSp A$ is the
hereditary coreflective hull of $\Kat A$.)
\end{PROP}

By $\FG$ we denote the category of all finitely generated spaces. It is well
known (see \eg \cite{HER}) that if $X$ is not finitely generated, then
$\FG\subseteq \CHb X$.

We say that a subcategory $\Kat A$ of $\Top$ is \emph{hereditary} if with each
\tsp $X$ it contains also all its \ssps. It is well known that the category of
all finitely generated spaces and all its subcategories that are coreflective
in $\Top$ are hereditary.

Some known hereditary coreflective subcategories of $\Top$ are $\Genal$ and
$\Topal$, where $\alpha$ is an infinite cardinal. $\Genal$ is the subcategory
of all spaces having tightness not exceeding $\alpha$. $\Topal$ is the category
of all \tsps such that the intersection of every family of open sets, which has
cardinality less than $\alpha$, is an open set.

Let $A$ be a \tsp. We say that $A$ is a \emph{\filsp{}} if it has precisely one
accumulation point. The following assertion is obvious.

\begin{LM}\label{LMSSPFSP}
Let $X$ be a \filsp with an accumulation point $a$ and let $Y$ be a \ssp of $X$
containing the point $a$, then the map $\Zobr fXY$, defined by $f(x)=x$ for
$x\in Y$ and $f(x)=a$ for $x\in X\setminus Y$, is a quotient map.
\end{LM}

Given a \tsp $X$ and a point $a\in X$, denote by $X_a$ the space constructed by
making each point, other than $a$, isolated with $a$ retaining its original
neighborhoods. (\Ie a subset $U\subseteq X$ is open in $X_a$ \iaoi $a\notin U$
or there exists an open subset $V$ of $X$ such that $a\in V\subseteq U$.) The
\tsp $X_a$ is called \emph{\primfac of $X$ at the point $a$.}  It is clear that
any \primfac is either a \filsp or a discrete space.

\begin{PROP}[{\cite[Proposition 3.5]{CIN}}]\label{PROPCIN}
If $\Kat A$ is a hereditary coreflective subcategory of $\Top$ with
$\FG\subseteq\Kat A$, then for each $X\in\Kat A$ and each $a\in X$ the \primfac
$X_a$ of $X$ at $a$ belongs to $\Kat A$.
\end{PROP}

Let $A$ be a \filsp with an accumulation point $a$. A \ssp $B$ of
$A$ is said to be a \emph{\filssp{}} of $A$ if $B$ is a \filsp
(\ie $a\in B$ and $\ol{B\setminus\{a\}}\ni a$).

\begin{LM}\label{LMNOVA2}
Let $(A_i; i\in I)$ be a family of \filsps and let $a_i\in A_i$ be an
accumulation point of $A_i$ for $i\in I$. A \tsp $X$ belongs to $\CHb{\{A_i;
i\in I\}}$ \iaoi for every non-closed subset $M$ of $X$ there exists $i\in I$,
a \filssp $B$ of $A_i$ and a continuous map $\Zobr fBX$ such that $\Obr
f{B\setminus\{a_i\}}\subseteq M$ and $f(a_i)\notin M$.
\end{LM}

\begin{proof}
Let $\Kat B\subseteq \Top$ be the class of all topological spaces satisfying
the given condition. First we show that $\Kat B$ is a coreflective subcategory
of $\Top$. It is evident that $\Kat B$ is closed under the formation of
topological sums. Now let $X\in\Kat B$ and $\Zobr qXY$ be a quotient map. Let
$M$ be a non-closed subset of $Y$. Then $\Invobr qM$ is a non-closed subset of
$X$, $X\in\Kat B$, so that there exists $i\in I$, a \filssp $B$ of $A_i$ and a
continuous map $\Zobr gBX$ such that $\Obr g{B\setminus\{a_i\}} \subseteq
\Invobr qM$ and $g(a_i)\notin \Invobr qM$. Then for $\Zobr{f=q\circ g}BX$ we
get $\Obr f{B\setminus\{a_i\}} \subseteq M$ and $f(a_i)\notin M$. Hence,
$Y\in\Kat B$ and $\Kat B$ is a coreflective subcategory of $\Top$.

Since evidently $A_i\in \Kat B$ for each $i\in I$, we have $\CHb{\{A_i; i\in
I\}}\subseteq \Kat B$. To prove the reverse inclusion we construct a quotient
map from a sum of subspaces of $A_i$ to arbitrary space $X\in \Kat B$. (Every
subspace of $A_i$ belongs to $\CHb{A_i}$ by Lemma \ref{LMSSPFSP}.)

Let $X\in\Kat B$. Let $\Zobr{f_j}{B_j}X$, $j\in J$, be the family of all
continuous maps such that $B_j$ is a \filssp of some $A_i$, $i\in I$. Let
$D(X)$ be the discrete space on the set $X$ and $\Zobr{id_X}{D(X)}X$ be the
identity map. It is easy to check that the map $\Zobr
f{D(X)\sqcup(\coprod_{j\in J}B_j)}X$ given by the maps $id_X$ and $f_j, j\in
J$, is a quotient map. 
\end{proof}

Cardinals are initial ordinals where each ordinal is the (well-ordered) set of
its predecessors. We denote the class of all ordinals by $\On$. If $\alpha$ is
a cardinal, then by $\alpha^+$ we denote the cardinal which is a successor of
$\alpha$.  A net in a \tsp defined on an ordinal $\alpha$ we call an
\emph{$\alpha$-sequence.}

From now on we assume that $A$ is a \filsp with an accumulation point $a$ which
is not finitely generated and the tightness of the space $A$ is $t(A)=\alpha$.

\section{Closure operator describing \CHA}

The notion of sequential closure was used in \cite{FRARAJ} when studying
sequential and subsequential spaces. Now we introduce a corresponding closure
operator for the subcategory $\CHA$.

Let $X$ be an arbitrary space and $M\subseteq X$. The set $M_1=\{x\in X:$ there
exists a \filssp $B$ of $A$ and a continuous map $\Zobr fBX$ such that $\Obr
f{B\setminus\{a\}} \subseteq M$ and $f(a)=x\}$ is called the \emph{$A$-closure
of $M$.} Using transfinite induction we can define the set $M_\beta$ (the
$\beta$-th $A$-closure of $M$) for each ordinal $\beta$ as follows. $M_0=M$,
$M_{\beta+1}=(M_\beta)_1$ for each ordinal $\beta$ and
$M_\gamma=\bigcup_{\beta<\gamma} M_\beta$ for each limit ordinal $\gamma>0$.
Put $\widetilde M=\bigcup_{\beta\in\On} M_\beta$.

Evidently $(\widetilde M)_1=\widetilde M$, $\widetilde M \subseteq \overline
M$. It is also clear that $M_\beta\subseteq M_\gamma$ holds for $\beta<\gamma$.
If $A\subseteq B\subseteq X$, then $A_\beta\subseteq B_\beta$ for each ordinal
$\beta$ and $\widetilde A\subseteq\widetilde B$. If $M_\beta=M_{\beta+1}$ for
some ordinal $\beta$, then $\widetilde M=M_\beta$.

The following proposition characterizes the spaces belonging to $\CHA$ using
the closure operator $M\mapsto \widetilde M$. It is a special case of
\cite[Theorem 3.1.7]{KAN} which includes more general cases of closure
operators.

\begin{PROP}\label{PROPTILDE}
A topological space $X$ belongs to $\CHA$ if and only if $\overline M =
\widetilde M$ for every subset $M\subseteq X$.
\end{PROP}

\begin{proof}
Let $X\in\CHA$ and $M\subseteq X$. Then $(\widetilde M)_1 \setminus \widetilde
M=\emptyset$, so that by Lemma \ref{LMNOVA2} $\widetilde M$ is closed and
$\widetilde M=\overline M$.

Conversely, if $\overline M=\widetilde M$ for each $M\subseteq X$ and $M$ is non-closed, then
$M_1\setminus M\neq\emptyset$ and there exists a \filssp $B$ of $A$ and a continuous map
$\Zobr fBX$ such that $\Obr f{B\setminus\{a\}} \subseteq M$ and $f(a)\notin M$. Hence,
according to Lemma \ref{LMNOVA2}, we conclude that $X\in\CHA$. 
\end{proof}

\begin{PROP} \label{PROPTIGH}
Let $A$ be a \filsp with an accumulation point $a$, $X\in\CHA$ and
$\alpha=t(A)$. Then for every subset $M\subseteq X$ it holds
$M_{\alpha^+}=\overline M$.
\end{PROP}

\begin{proof}
If suffices to prove that $(M_{\alpha^+})_1=M_{\alpha^+}$. Let $c\in
(M_{\alpha^+})_1$. Then there exists a \filssp $B$ of $A$ and a continuous map
$\Zobr fBX$ with $f(a)=c$ and $\Obr f{B\setminus\{a\}}\subseteq M_{\alpha^+}$.
Since $t(A)=\alpha$ and $a\in\ol{B\setminus\{a\}}$, there exists $C\subseteq
B\setminus\{a\}$ with $\card C\leq\alpha$ such that $a\in\ol C$. The \ssp
$B_1=C\cup\{a\}$ of $A$ is a \filssp, $\Zobr{f|_{B_1}}{B_1}X$ is continuous and
$\Obr{f|_{B_1}}C\subseteq M_{\alpha^+}$.

For each $x\in C$ choose $\beta_x < \alpha^+$ such that $x\in M_{\beta_x}$
($\alpha^+$ is a limit ordinal). Since $\card C\leq\alpha<\alpha^+$ and
$\alpha^+$ is \arcn we obtain that $\gamma=\sup\{\beta_x, x\in C\}<\alpha^+$.
Then $C\subseteq M_\gamma$ and, obviously, $f|_{B_1}(a)=f(a)=c \in M_{\gamma+1}
\subseteq M_{\alpha^+}$. Thus, $(M_{\alpha^+})_1 \subseteq M_{\alpha^+}$. 
\end{proof}

\section{\filsum A{}}

The notion of \filsum A is a special case of the brush defined in \cite{KAN}
and a generalization of the sequential sum introduced in \cite{ARHFRA}. The
sequential sum was used in \cite{FRARAJ} for constructing the set of
\uv{canonical} \filsps which generates the category of subsequential spaces.
The notion of the \filsum A will be used in a similar way to produce the set of
special \filsps that generates $\SCHA$.

\begin{DEF}
Let $A$ be a \filsp with an accumulation point $a\in A$. Let us denote
$B:=A\setminus\{a\}$. Let for each $b\in B$ $X_b$ be a topological space and
$x_b\in X_b$. Then the \emph{\filsum A} $\SumA A{X_b}{x_b}$ is the topological
space on the set $F=A\cup(\bigcup\limits_{b\in B} \{b\}\times (X_b\setminus
\{x_b\}))$ such that the map $\Zobr{\varphi}{A\sqcup(\coprod\limits_{b\in B}
X_b)}F$ given by $\varphi(x)=x$ for $x\in A$, $\varphi(x)=(b,x)$ for $x\in
X_b\setminus\{x_b\}$ and $\varphi(x_b)=b$ for every $b\in B$ is a quotient map.
(We assume $A$ and all $\{b\}\times X_b$ to be disjoint.) The map $\varphi$
will be called the \emph{defining map} of the \filsum A.
\end{DEF}

Often it will be clear from the context what we mean under $A$ and we will
abbreviate the notation of the \filsum A to $\Suma{X_b}{x_b}$ or $\Sum{X_b}$.
The \filsum A is obtained simply by identifying every $x_b\in X_b$ with the
point $b\in A$. It is easy to see that the subspace $\Obr \varphi{X_b}$ is
homeomorphic to $X_b$ and $A$ is also a subspace of the \filsum A
$\Suma{X_b}{x_b}$.

The \filsum A is defined using topological sum and quotient map, thus if $\Kat
A$ is a coreflective subcategory of $\Top$ and $\Kat A$ contains $A$ and all
$X_b$'s, then the \filsum A $\Sum X_b$ belongs to $\Kat A$.

The following lemma follows easily from the definition of the \filsum A.

\begin{LM}\label{LMOPACLINSUM}
A subset $U\subseteq$ {\SumA A{X_b}{x_b}} is open (closed) if and only if
$U\cap A$ is open (closed) in $A$ and $U\cap\Obr\vp{X_b}$ is open (closed) in
$\Obr\vp{X_b}$ for every $b\in B$.
\end{LM}

Let for every $b\in B$ $X_b$ and $Y_b$ be topological spaces, $x_b\in X_b$,
$y_b\in Y_b$ and let $\Zobr{f_b}{X_b}{Y_b}$ be a function with $f_b(x_b)=y_b$.
Then we can define a map $\Zobr{f=:\Sum{f_b}}{\SumA A{X_b}{x_b}}{\SumA
A{Y_b}{y_b}}$ by $f(b,x)=(b,f_b(x))$ for $x\in X_b\setminus \{x_b\}$ and
$f(x)=x$ for $x\in A$. Let us note that $f \circ \varphi_1|_{X_b} =
\varphi_2|_{Y_b} \circ f_b$ where $\varphi_1$ and $\varphi_2$ are the defining
maps of the \filsums A $\Sum{X_b}$ and $\Sum{Y_b}$ respectively.

We will need the following simple lemma:

\begin{LM} \label{LMPOM}
Let $\Zobr fXY$ be a quotient map, $A\subseteq Y$ and let $f$ be one-to-one
outside $A$. Then $\Zobr{f|_{\Invobr fA}}{\Invobr fA}A$ is a quotient map.
\end{LM}

\begin{LM} \label{LMSUMAAEMB}
Let $A$ be a \filsp with an accumulation point $a$ and $B=A\setminus\{a\}$. Let
for every $b\in B$ $\Zobr{f_b}{X_b}{Y_b}$ be a map between topological spaces,
$x_b\in X_b$, $y_b\in Y_b$ and $f_b(x_b)=y_b$.
\begin{enumerate}
\renewcommand{\theenumi}{\roman{enumi}}
\renewcommand{\labelenumi}{(\theenumi)}
 \item If all $f_b$'s are continuous, then $\Sum{f_b}$ is continuous.
 \item If all $f_b$'s are quotient maps, then $\Sum{f_b}$ is a quotient map.
 \item If all $f_b$'s are embeddings, then $\Sum{f_b}$ is an embedding.
 \item If all $f_b$'s are homeomorphisms, then $\Sum{f_b}$ is a homeomorphism.
 \item Let $C$ be a \filssp of $A$. Then $\SumA C{X_b}{x_b}$ is a \ssp
   of the space $\SumA A{X_b}{x_b}$.
\end{enumerate}
\end{LM}

\begin{proof}
Put $f=\Sum{f_b}$ and let $\varphi_1$, $\varphi_2$ be the defining maps of the
\filsums A $\Suma{X_b}{x_b}$, $\Suma{Y_b}{y_b}$ respectively. Let us denote
$id_A\sqcup(\coprod_{b\in B} f_b)$ by $h$. In this situation the following
diagram commutes.
$$\xymatrix{
  A\sqcup(\coprod X_b) \ar[r]^{h} \ar[d]_{\varphi_1} &
    A\sqcup(\coprod Y_b \ar[d]_{\varphi_2}) \\
  {\Suma{X_b}{x_b}} \ar[r]_f & \Suma{Y_b}{y_b}
}$$ The validity of (i) and (ii) follows easily from the fact that $\varphi_1$
and $\varphi_2$ are quotient maps.

(iii) Now, suppose that all $f_b$'s are embeddings. \Wlog we can assume that
$X_b\subseteq Y_b$ and $f_b$ is the inclusion of $X_b$ into $Y_b$ for every
$b\in B$. Let $X'$ be the \ssp of the space $\Sum{Y_b}$ on the set $\Sum{X_b}$.
We have the following situation:
$$\xymatrix{
  A\sqcup(\coprod X_b) \ar@{^{(}->}[r]^{h} \ar[d]_{\varphi_1} & A\sqcup(\coprod Y_b \ar[d]_{\varphi_2}) \\
  X' \ar@{^{(}->}[r]_f & \Sum{Y_b}
}$$ We only need to prove that $X'$ has the quotient topology with respect to
$\varphi_1$, because this implies that $X'=\Sum{X_b}$ and $f$ is an embedding
of $X'=\Sum{X_b}$ to $\Sum{Y_b}$. But $\varphi_2$ is one-to-one outside the set
$A \sqcup (\coprod X_b)$ and Lemma \ref{LMPOM} implies that $\varphi_1$ is a
quotient map.

(iv) It is an easy consequence of (ii) and (iii). (v) It follows easily from
the definition of the \filsum A. 
\end{proof}

\begin{COR} \label{CRSSPSUM}
Let $A$ be a \filsp with an accumulation point $a$ and let $C$ be a \filssp of
$A$. Let for every $b\in A\setminus \{a\}$ $X_b$ be a \tsp and $x_b\in X_b$.
Let for every $b\in C$ $Y_b$ be a \ssp of $X_b$ such that $x_b\in Y_b$. Then
$\SumA C{Y_b}{x_b}$ is a \ssp of the space $\SumA A{X_b}{x_b}$.
\end{COR}

Let us note, that if for every $b\in A\setminus\{a\}$ $f_b$ is an embedding
which maps isolated points of $X_b$ to isolated points of $Y_b$, then the
embedding $\Sum{f_b}$ has the same property.

\section{The sets $TS_\gamma$, $TSS_\gamma$}

In this section we construct the set of special \filsps that generates $\SCHA$
(where $A$ is a \filsp which is not finitely generated and $t(A)=\alpha$). We
start with defining the set $TS_\gamma$ of \tsps for each ordinal
$\gamma<\alpha^+$.

Let $TS_0=\emptyset$ and $TS_1$ be the set of all \filssps of $A$.

If $\beta\geq1$ is an ordinal, then $TS_{\beta+1}$ consists of all \filsums B
$\SumA B{X_b}{x_b}$ where $B$ is a \filssp of $A$, each $X_b\in TS_\beta$ and
$x_b=a$.

If $\gamma>0$ is a limit ordinal, then $TS_\gamma=\bigcup_{\beta<\gamma}
TS_\beta$.

Sometimes, if we want to emphasize which \filsp $A$ is used to construct this
set, we use the notation $TS_\gamma(A)$.

Every space belonging to $TS_\gamma$ contains $B$ as a \ssp and therefore it
contains $a$. All spaces from $TS_\gamma$ are constructed from $A$ using
\filsums B, where $B\in\CHA$, thus $TS_\gamma \subseteq \CHA$ for each
$\gamma$.

The following lemma is a generalization of \cite[Lemma 6.2]{FRARAJ}.

\begin{LM}\label{LMTSCLOP}
Let $X$ be a \tsp and $M\subseteq X$. If $p\in M_\beta\setminus M_\gamma$ for
any $\gamma<\beta$, then there exists a space $S\in TS_\beta$ and a continuous
map $\Zobr f{S}X$, which maps all isolated points of $S$ into $M$ and maps only
the point $a$ to $p$.
\end{LM}

\begin{proof}
For $\beta=1$ the claim follows from the definition of $M_1$.

From the definition of $M_\beta$ it follows that $\beta$ is a non-limit
ordinal. According to Proposition \ref{PROPTIGH} $\beta<\alpha^+$. Suppose the
assertion is true for any subset $K$ of $X$ and for any $\beta'<\beta$.

For a non-limit $\beta>1$ there exists a \filssp $B$ of $A$ and a continuous
map $\Zobr fBX$ such that $f(a)=p$ and $\Obr f{B\setminus\{a\}} \subseteq
M_{\beta-1}$.

If $\beta-1$ is non-limit, we can moreover assume that $\Obr f{B\setminus\{a\}}
\subseteq M_{\beta-1}\setminus M_{\beta-2}$. (If necessary, we choose
$B'=\{b\in B: f(b)\in M_{\beta-1}\setminus M_{\beta-2}\}$ and $f'=f|_{B'}$.
$B'$ is a \filssp of $A$, otherwise we get $x\in M_{\beta-1}$.)

If $\beta-1$ is a limit ordinal, then for each point $x\in M_{\beta-1}$ there
exists the smallest ordinal $\gamma<\beta-1$ such that $x\in M_\gamma$.
Obviously, $\gamma$ is a non-limit ordinal.

Thus for each $x\in \Obr f{B\setminus\{a\}}$ there exists a continuous map
$\Zobr{f_x}{S_x}X$, where $S_x\in TS_{\beta-1}$, which sends all isolated
points of $S_x$ into $M$ and $a$ to $x$.

Then $\SumA B{S_{f(b)}}a \in TS_\beta$ and we can define a map $\Zobr g{\SumA
B{S_{f(b)}}a}X$ such that $g|_B=f$ and $g|_{\{x\}\times
(S_x\setminus\{a\})}(x,y)=f_x(y)$ for $y\in S_x\setminus\{a\}$. Clearly, $g$
maps isolated points into $M$. It remains only to show that $g$ is continuous.

The defining map $\Zobr \vp{B\sqcup(\coprod_{b\in B\setminus\{a\}}S_{f(b)})}{\Suma
{S_{f(b)}}a}$ is a quotient map. Thus, $\Zobr g{\Suma {S_{f(b)}}a}X$ is continuous \iaoi
$g\circ\vp$ is continuous. But $g\circ\vp|_{B}=f$ and $g\circ\vp|_{S_x}=f_x$ are continuous,
thus $g$ is continuous. 
\end{proof}

For any $S \in TS_\gamma$ we denote by $P(S)$ the \ssp of the space $S$ which
consists of all isolated points of $S$ and of the point $a$. Clearly, $P(S)$ is
a \filsp. We denote by $TSS_\gamma$ the set of all spaces $P(S)$ where $S\in
TS_\gamma$. The above lemma implies:

\begin{LM} \label{LMTSSCLOP}
If $p\in M_\beta$ and $p\notin M_\gamma$ for any $\gamma<\beta$, then there
exists a space $T\in TSS_\beta$ and a continuous map $\Zobr f{T}X$, which maps
all isolated points of the space $T$ into $M$ and such that $f(a)=p$.
\end{LM}

\begin{PROP}\label{PRGENSCHATH}
$\SCHA$ is generated by the set $\bigcup\limits_{\gamma<\alpha^+} TSS_\gamma$.
\end{PROP}

\begin{proof}
Let $X\in\SCHA$. According to Lemma \ref{LMNOVA2} it suffices to prove that for
any subset $M\subseteq X$ and any $x\in\ol M\setminus M$ there exists
$T\in\bigcup\limits_{\gamma<\alpha^+} TSS_\gamma$ and a continuous map $\Zobr
fTX$ such that $f(a)=x$ and $\Obr f{T\setminus\{a\}}\subseteq M$.

Since $X\in\SCHA$ there exists $Y\in\CHA$ such that $X$ is a \ssp of $Y$. Denote by $\Uzav
MY$ the closure of $M$ in $Y$. Then $\ol M=\Uzav MY \cap X$ and $x\in \Uzav MY\setminus M$ in
$Y$. By Proposition \ref{PROPTIGH} $\Uzav MY=M_{\alpha^+} = \bigcup\limits_{\beta<\alpha^+}
M_\beta$. Let $\beta$ be the smallest ordinal with $x\in M_\beta$. Then $\beta>0$ and for any
$\gamma<\beta$ $x\notin M_\gamma$. By Lemma \ref{LMTSCLOP} there exists $S\in TS_\beta$ and a
continuous map $\Zobr fSY$ with $f(a)=x$ and $f(c)\in M$ for any isolated point of $S$. Then
$P(S)\in TSS_\gamma$ and $\Obr f{P(S)}\subseteq X$. Hence, $\Zobr {f|_{P(S)}}{P(S)}X$ is a
continuous map satisfying the required conditions. Consequently,
$X\in\CHb{\bigcup\limits_{\gamma<\alpha^+} TSS_\gamma}$. 
\end{proof}

\begin{REM}
It can be easily seen that if we define the sets $T'S_\gamma$,
$\gamma<\alpha^+$, similarly as the sets $TS_\gamma$ but we use only the
\filsums A (and not all \filsums B for \filssps $B$ of $A$) and then we put
$T'SS_\gamma=\{P(S): S\in T'S_\gamma\}$ we obtain the set
$\bigcup\limits_{\gamma<\alpha^+} T'SS_\gamma$ which also generates $\SCHA$.
This follows from the fact that any space from
$\bigcup\limits_{\gamma<\alpha^+} TSS_\gamma$ is a \filssp of some space from
$\bigcup\limits_{\gamma<\alpha^+} T'SS_\gamma$.

Similarly, if we put $\ciTSS\gamma=\{S_a: S\in T'S_\gamma\}$ ($S_a$ is the
\primfac of $S$ at $a$), then the set
$\bigcup\limits_{\gamma<\alpha^+}\ciTSS\gamma$ generates $\SCHA$ because
$\bigcup\limits_{\gamma<\alpha^+}\ciTSS\gamma \subseteq \SCHA$ and for every
$S\in\bigcup\limits_{\gamma<\alpha^+} T'S_\gamma$ $P(S)$ is a \ssp of $S_a$.
\end{REM}

\section{The spaces $\Som$ and $(\Som)_a$}

The space $\Som$ is defined similarly as $S_\omega$ in \cite{ARHFRA} using the
\filsum A and the space $A$ instead of  the sequential sum and the space
$\Com$. We start with defining the space $A_n$ for each $n\in\N$ putting
$A_1=A$ and $A_{n+1}=\SumA A{A_n}a$. Clearly, $A_1$ is a \ssp of $A_2$ and if
$A_{n-1}$ is a \ssp of $A_n$, then, according to Lemma \ref{LMSUMAAEMB},
$A_n=\SumA A{A_{n-1}}a$ is a \ssp of $A_{n+1}=\SumA A{A_n}a$. Hence, $A_n$ is a
\ssp of $A_{n+1}$ for each $n\in\N$.

\begin{figure}[h]
\centerline{\epsfig{file=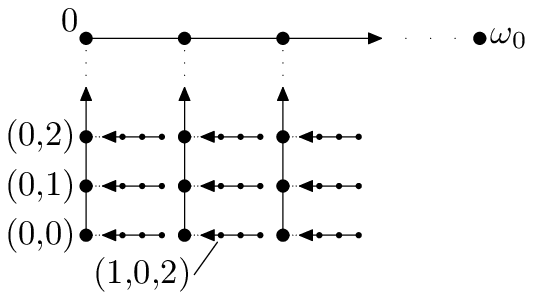}}
\caption{The space $A_3$ for $A=\Com$}\label{FIG1}
\end{figure}

The Figure \ref{FIG1} represents the space $A_3$ for $A=\Com$.
(The space $\Com$ is defined in Example \ref{EXACOM}.)

The space $\Som$ is a \tsp defined on the set $\bigcup_{n\in\N} A_n$ such that
a subset $U$ of $\bigcup_{n\in\N} A_n$ is open in $\Som$ \iaoi $U\cap A_n$ is
open in $A_n$ for every $n\in\N$. It is obvious that for every $n\in\N$ the
space $A_n$ is a \ssp of $\Som$ and $\Som$ is a quotient space of the
topological sum $\coprod_{n\in\N} A_n$. Consequently, $\Som$ belongs to $\CHA$.
Observe that $\Som$ can be considered as an inductive limit of its \ssps $A_n$,
$n\in\N$.

Similarly as the space $S_\omega$ in \cite{ARHFRA} the space $\Som$ has the
following important property.

\begin{PROP}
$\Som=\SumA A{\Som}a$
\end{PROP}

\begin{proof}
Put $X=\SumA A{\Som}a$. For each $n\in A$ the space $A_n$ is a \ssp of $\Som$
and it follows that $A_{n+1}=\SumA A{A_n}a$ is a \ssp of $X$ (Lemma
\ref{LMSUMAAEMB}). Obviously, $A=A_1$ is also a \ssp of $X$ and we obtain that
for each $n\in\N$ $A_n$ is a \ssp of $X$. Clearly, $X=\bigcup\limits_{n\in\N}
A_n$. To finish the proof it suffices to check that if $U$ is a subset of $X$
and $U\cap A_n$ is open in $A_n$ for each $n\in\N$, then $U$ is open in
$X$.

Let us denote by $A_n^b$ the \ssp of $X$ on the set
$\{b\}\cup(\{b\}\times(A_n\setminus\{\alpha\}))$ and by $\Som^b$ the \ssp of
$X$ on the set $\{b\}\cup(\{b\}\times(\Som\setminus\{\alpha\}))$. Clearly,
$A_n^b$ is homeomorphic to $A_n$ and $\Som^b$ is homeomorphic to $\Som$,
$A_n^b$ is a \ssp of $\Som^b$ and a subset $V$ of $\Som^b$ is open in $\Som^b$
\iaoi $V\cap A_n^b$ is open in $A_n^b$ for each $n\in\N$.

If $U\subseteq X$ and $U\cap A_n$ is open in $A_n$ for all $n\in\N$, then
$U\cap A$ is open in $A$ and $U\cap A_{n+1}$ is open in $A_{n+1}=\SumA A{A_n}a$
for all $n\in\N$. Then $U\cap A_n^b$ is open in $A_n^b$ for each $n\in\N$ and
$b\in A\setminus\{a\}$ and  it follows that $U\cap\Som^b$ is open in $\Som^b$
for each $b\in A\setminus \{a\}$. Hence, $U$ is open in $X$. 
\end{proof}

The following lemma is evident.

\begin{LM}
$\card{\Som}=\card A$
\end{LM}

\begin{LM}\label{LMSSPSOM}
For every ordinal $\gamma$, $1\leq\gamma<\alpha^+$ and every space $S\in
TS_\gamma$ the space $S$ is a subspace of $\Som$. (Clearly, the point $a$ of
$S$ coincides with the point $a$ of $\Som$.)
\end{LM}

\begin{proof}
If $\gamma=1$, then $S=B$ is a \filssp of $A=A_1$. Let $\gamma$ be an ordinal,
$1<\gamma<\alpha^+$ and the assertion hold for every ordinal $\beta$,
$1\leq\beta<\gamma$. If $S=\SuM B{X_b} \in TS_\gamma$, then for each $b\in
B\setminus\{a\}$ $X_b \in TS_{\beta_b}$ with $1\leq\beta_b<\gamma$. Hence, for
each $b\in B\setminus\{a\}$, $X_b$ is a \ssp of $\Som$  and, according to
Corollary \ref{CRSSPSUM}, $S$ is a \ssp of $\Som=\SumA A{\Som}a$. 
\end{proof}

\begin{THM}\label{THM54}
Let $(\Som)_a$ be the \primfac of the space $\Som$ at $a$. Then $(\Som)_a$ is a
\filsp, $\CHb{(\Som)_a}=\SCHA$ and $\card (\Som)_a=\card A$.
\end{THM}

\begin{proof}
Evidently, $(\Som)_a$ is a \filsp and $\card (\Som)_a=\card A$.
Since $\Som$ belongs to $\CHA\subseteq\SCHA$, according to
Proposition \ref{PROPCIN} $(\Som)_a$ belongs to $\SCHA$. Hence, it
suffices to check that $\bigcup\limits_{\gamma<\alpha^+}
TSS_\gamma \subseteq \CHb{(\Som)_a}$.

Let $T\in \bigcup\limits_{\gamma<\alpha^+} TSS_\gamma$. Then there exists an
ordinal $\gamma$, $1\leq\gamma<\alpha^+$, and $S\in TS_\gamma$ such that
$T=P(S)$. By Lemma \ref{LMSSPSOM} $S$ is a \ssp of $\Som$ and, clearly, it
follows that $T=P(S)$ is a \ssp of $(\Som)_a$. Consequently, there exists a
quotient map $(\Som)_a \to T$ and we obtain that $T$ belongs to
$\CHb{(\Som)_a}$. 
\end{proof}

Finally, let $X$ be an arbitrary \tsp which is not finitely generated and
$\{X_c, c\in Y\}$ be the set of all \primfacs of $X$ that are not discrete
spaces. Denote by $A_X$ the quotient space of the topological sum
$\coprod\limits_{c\in Y}(\{c\}\times X_c)$ obtained by collapsing all points of
the subset $\{(c,c), c\in Y\}$ of the space $\coprod\limits_{c\in
Y}(\{c\}\times X_c)$ into one point $a$.

The space $A_X$ is a \filsp which is not finitely generated, $a$ is the
accumulation point of $A_X$, $\card A_X=\card X$ and the following statement
holds:

\begin{THM}\label{THMSIMPLEGEN}
$\SCHb X=\CHb{((A_X)_\omega)_a}$, $((A_X)_\omega)_a$ is a \filsp
and $\card ((A_X)_\omega)_a=\card X$.
\end{THM}

\begin{proof}
Evidently, $A_X\in\SCHb X$, $X\in\CHb{A_X}$ and therefore  $\SCHb
X=\SCHb{A_X}$. According to Theorem \ref{THM54}
$\SCHb{A_X}=\CHb{((A_X)_\omega)_a}$, $((A_X)_\omega)_a$ is a \filsp and $\card
((A_X)_\omega)_a=\card A_X=\card X$. 
\end{proof}

Recall, that a \tsp $X$ belongs to $\Topom$ \iaoi every countable intersection
of open subsets of $X$ is open in $X$ and $\Topom$ is a hereditary coreflective
subcategory of $\Top$. If the space $A$ belongs to $\Topom$, then $\SCHA
\subseteq \Topom$ and we can find smaller (and simpler) set of generators of
$\SCHA$ than the set $\bigcup_{\gamma<\alpha^+} TSS_\gamma(A)$ constructed in
Proposition \ref{PRGENSCHATH}.

\begin{PROP}\label{PROPAMNAOM}
If $A\in\Topom$, then $\SCHA=\CHb{\{\Anm; 0<n<\omnul\}}$.
\end{PROP}

\begin{proof}
It suffices to show that $(\Aom)_a \in \CHb{\{\Anm; n<\omnul\}}$. Clearly, each
$\Anm$ is a \ssp of $(\Aom)_a$. Denote by $\Vloz{i_n}{\Anm}{(\Aom)_a}$ the
corresponding embedding and by
$\Zobr{f}{\coprod\limits_{n\in\N}\Anm}{(\Aom)_a}$ the continuous map given by
the maps $i_n$, $n\in\N$. It is easy to see that this map is surjective. We
claim that $f$ is also a quotient map.

It suffices to show that if $a\in U\subseteq (\Aom)_a$ and $U\cap\Anm$ is open
in $\Anm$ for each $n$, $0<n<\omnul$, then $U$ is open in $(\Aom)_a$. Since
$\Anm$ is a \ssp of $\Aom$, there exists an open subset $W_n$ of $\Aom$ such
that $W_n\cap\Anm=U\cap\Anm$. Put $W=\bigcap_{0<n<\omnul} W_n$. The set $W$ is
open in $\Aom$ since $\Aom$ belongs to $\Topom$. We have $W\cap\Anm \subseteq
U\cap\Anm$ and $\bigcup_{0<n<\omnul}\Anm=\Aom$, therefore $W\subseteq U$.
Obviously, $a\in W$. Hence, the set $U$ is open in $(\Aom)_a$.
\end{proof}

Next we present some special cases of our construction.

\begin{EXA}\label{EXACOM}
\emph{Sequential spaces.} Recall that subspaces of sequential spaces are called
subsequential. The category $\Seq$ of sequential spaces is the coreflective
hull of the space $\Com$. The space $\Com$ is the \tsp on the set
$\omega_0+1=\omega_0\cup\{\omega_0\}$ such that all points of $\omega_0$ are
isolated and a set containing $\omega_0$ is open \iaoi its complement is
finite. (Equivalently, the topology of $\Com$ is the order topology given by
the usual well-ordering of $\omega_0+1$.) The space $\Com_\omega$  is
homeomorphic to $S_\omega$ defined in \cite{ARHFRA}. Our results imply that the
\primfac of the space $\Com_\omega$ at $\omega_0$ is a generator of the
category of subsequential spaces. Another countable generator of this category
was constructed before in \cite{FRARAJ}.
\end{EXA}

\begin{EXA}\label{EXACALP}
\emph{The coreflective hull of the space $\Calp$.} Let $\alpha$ be \arcn and
$\Calp$ be the \tsp on the set $\alpha+1=\alpha\cup\{\alpha\}$ such that all
points of $\alpha$ are isolated and a set containing $\alpha$ is open \iaoi its
complement has cardinality less than $\alpha$. It is well known that $X$
belongs to $\CHb{\Calp}$ \iaoi a subset $V\subseteq X$ is closed in $X$
whenever for each $\alpha$-sequence of points from $V$ the set $V$ contains
also all limits of this $\alpha$-sequence. The subcategories $\SCHb{\Calp}$ are
minimal elements of the collection of all hereditary coreflective subcategories
of $\Top$ above $\FG$. We use the subcategories $\SCHb{\Calp}$ in the next
section. Our construction yields the generator $(\Calp_\omega)_\alpha$ of
$\SCHb{\Calp}$ which has cardinality $\alpha$.
\end{EXA}

\section{Subcategories of $\Top$ having $\FG$ as their hereditary coreflective kernel}

Recall that a \emph{hereditary coreflective kernel} of a subcategory $\Kat A$
of $\Top$ is the largest hereditary coreflective subcategory of $\Top$
contained in $\Kat A$. We denote it by $\HCK{A}$. In this section we prove that
if $\Kat A$ and $\Kat B$ are coreflective subcategories of $\Top$ such that
$\HCK A=\HCK B=\FG$, then also $\HCKb{\CHb{\Kat A\cup \Kat B}}=\FG$. The
analogous result does not hold for infinite countable joins of coreflective
subcategories. This problem is closely related to the subcategories
$\SCHb{\Calp}$ because (see \cite[Theorem 4.8]{CIN}) $\FG$ is the hereditary
coreflective kernel of a coreflective subcategory $\Kat A$ of $\Top$ \iaoi
$\FG\subseteq\Kat A$ and for any \rcn $\alpha$ the category $\SCHb{\Calp}$ is
not contained in $\Kat A$.

In \cite[Problem 7]{HERHUS} H.~Herrlich and M.~Hu\v{s}ek suggest to study
coreflective subcategories of $\Top$ such that their hereditary coreflective
hull is the whole category $\Top$ (\ie $\SSp A=\Top$) and their hereditary
coreflective kernel is the subcategory $\FG$. In the paper \cite{SLEZ} it is
shown that there exists the smallest such subcategory of $\Top$ and the
collection of all such subcategories of $\Top$ is closed under the formation of
arbitrary non-empty intersections. In this section we prove that this
collection is also closed under the formation of non-empty finite joins without
being closed under the formation of infinite countable joins in the lattice of
all coreflective subcategories of \Top.

Throughout this section we will apply the results obtained in preceding
sections to \filsps $\Calp$, $\alpha$ being \arcn, defined in Example
\ref{EXACALP}. Note that $\alpha$ is an accumulation point of $\Calp$ and
$t(\Calp)=\alpha$ for any \rcn $\alpha$. Since any \filssp of $\Calp$ is
homeomorphic to $\Calp$ it suffices to use only $\Calp$-sums in the definition
of $TS_\gamma$. For instance, if $n$ is a natural number, then $TS_n$ as well
as $TSS_n$ contain precisely one space.

In order to prove the main result of this section, we first prove that if
$\SCHb{\Calp}\subseteq \CHb{\Kat A\cup \Kat B}$ for some coreflective subcategories $\Kat A$,
$\Kat B$ of $\Top$, then one of these subcategories contains $\SCHb{\Calp}$. We show it
separately for the case $\alpha=\omnul$ and $\alpha\geq\omega_1$.

We start with the case $\alpha=\omnul$ where we can use some results presented
in the paper \cite{FRARAJ}. As the sets $TSS_\gamma$, $\gamma<\omega_1$,
introduced in \cite{FRARAJ} do not coincide with the sets $TSS_\gamma(\Com)$
defined in Section 4 we denote the sets used in \cite{FRARAJ} by $\TSSd\gamma$.

The next lemma follows from \cite[Theorem 7.1]{FRARAJ}, \resp \cite[Corollary 7.2]{FRARAJ}.

\begin{LM}\label{LMFRA}
The category $\SSeq=\SCHb{\Com}$ of subsequential spaces is the
coreflective hull of the set $\bigcup_{\gamma<\omega_1}
\TSSd\gamma$.
\end{LM}

As a consequence of \cite[Theorem 7.1]{FRARAJ} and \cite[Theorem 6.4]{FRARAJ}
we obtain:

\begin{LM}\label{LMFRB}
If $\beta<\gamma<\omega_1$, then $\TSSd\beta \subseteq \CHb{\TSSd\gamma}$.
\end{LM}

The following result concludes the part of this section concerning the
subcategory $\SCHb{\Com}$.

\begin{PROP}\label{PROPCOMI0}
If $\SCHb{\Com} \subseteq \CHb{\bigcup_{i\in I} \Kat A_i}$, $\Kat A_i$ is a
coreflective subcategory of $\Top$ for every $i\in I$ and $\card I\leq \omnul$,
then there exists $i_0\in I$ such that $\SCHb\Com\subseteq \Kat A_{i_0}$.
\end{PROP}

\begin{proof}
Put $\beta_i=\sup\{\beta: \TSSd\beta\subseteq \Kat A_i\}$ for $i\in I$. Since $\sup
\beta_i=\omega_1$ (Lemma \ref{LMFRA}) and $\omega_1$ is \arcn, there exists $i_0\in I$ such
that $\beta_{i_0}=\omega_1$. By Lemma \ref{LMFRA} and Lemma \ref{LMFRB} we get that the
coreflective subcategory $\Kat A_{i_0}$ contains the subcategory $\SSeq=\SCHb\Com$.
\end{proof}

Next we want to prove a result analogous to Proposition \ref{PROPCOMI0} for the
space $\Calp$, where $\alpha\geq\omega_1$ is \arcn. In the case
$\alpha\geq\omega_1$ the desired result holds only for non-empty finite joins
of coreflective subcategories of $\Top$.

Recall that $\Calp_1=\Calp$ and $\Calp_{n+1}=\SumA \Calp{\Can}\alpha$.
According to Corollary \ref{CRSSPSUM} we obtain that $P(\Calp_{n+1})=P(\SumA
\Calp{P(\Can)}\alpha)$ and it is easy to see that $\alpha^{n+1}\cup\{\alpha\}$
is the underlying set of the space $P(\Calp_{n+1})$ and the \ssp of $\sum
P(\Can)$ on the set $\{\eta\} \cup (\{\eta\} \times \alpha^n)$ is homeomorphic
to $P(\Can)$ for each $\eta<\alpha$. To simplify the notation we will write
$\Canm$ instead of $P(\Can)$.

The following result is a special case of Proposition \ref{PROPAMNAOM}.

\begin{PROP}\label{PROPAMNAOMCALP}
If $\alpha\geq\omega_1$ is \arcn, then $\SCHb{\Calp}=\CHb{\{\Canm;
0<n<\omnul\}}$.
\end{PROP}

\begin{LM}\label{LMSA}
Let $\alpha\geq\omega_1$ be \arcn. If $M$ is a subset of $\Can$
such that $\alpha\in\ol M$ and $M$ contains only isolated points
of $\Can$, then there exists a subset $M'\subseteq M$ such that
the \ssp of the space $\Can$ on the set $\ol{M'}$ is homeomorphic
to $\Can$.
\end{LM}

\begin{proof}
The case $n=1$ is clear. Let the assertion be true for $m$. Denote the \ssp of
$\Camj=\SuM{\Calp}{\Cam}$ on the set $\{\eta\} \cup (\{\eta\}\times
(\Cam\setminus\{\alpha\}))$, where $\eta<\alpha$, by $\Camk$.

Put $B=\ol M \cap \Calp$. Then $B$ is a \filssp of $\Calp$, for each $\eta\in
B\setminus\{\alpha\}$ all points of the set $M_\eta=M\cap \Camk$ are isolated
in the space $\Camk$ and $\eta\in\ol{M_\eta}$ in $\Camk$ (observe that
$\ol{M_\eta}$ in $\Camk$ coincides with $\ol{M_\eta}$ in $\Camj$ because
$\Camk$ is closed in $\Camj$). Since $\Camk$ is homeomorphic to $\Cam$ by the
induction assumption we obtain that there exists a subset $M'_\eta\subseteq
M_\eta$ such that $\eta\in\ol{M'_\eta}$ and the \ssp $\ol{M'_\eta}$ of $\Camk$
is homeomorphic to some space $\Cam$.

Let $B'=B\setminus\{\alpha\}$ and $M'=\bigcup\limits_{n\in B'} M'_\eta$.
Clearly, $M'\subseteq M$, $\ol{M'}=\bigcup\limits_{\eta\in B'} \ol{M'_\eta}
\cup \{\alpha\}$ in $S$ and $\ol{M'_\eta}$ is homeomorphic to $\Cam$ for each
$\eta\in B'$.

The \ssp $B$ of $\Calp$ is homeomorphic to $\Calp$ and it is easy to check that
$\ol{M'}$ is homeomorphic to $\SuM\Calp\Cam=\Camj$. 
\end{proof}

\begin{COR}\label{CORSA}
Let $\alpha\geq\omega_1$ be \arcn, $0<n<\omnul$. Then every \filssp $T$ of
$\Canm$ is homeomorphic to $\Canm$.
\end{COR}

\begin{proof}
Put $M=T\setminus\{\alpha\}$. Clearly, $\alpha\in\ol M$. According to Lemma
\ref{LMSA} there exists a subset $M'$ of $M$ such that the \ssp
$M'\cup\{\alpha\}$ of $\Canm$ is homeomorphic to $\Canm$. It follows from the
proof of Lemma \ref{LMSA} that $M\setminus M'$ is a discrete clopen \ssp of
$\Canm$ with cardinality at most $\alpha$. Hence, $T=M\cup\{\alpha\}$ is
homeomorphic to $\Canm$ as well. 
\end{proof}

\begin{PROP}\label{PRAI0}
Let $\alpha\geq\omega_1$ be \arcn and $0<n<\omnul$. If $\Canm\in
\CHb{\bigcup_{i\in I} \Kat A_i}$, where all $\Kat A_i$'s are coreflective
subcategories of $\Top$, then there exists $i_0\in I$ such that $\Canm\in\Kat
A_{i_0}$.
\end{PROP}

\begin{proof}
The space $\Canm$ is a \filsp with an accumulation point $\alpha$. If $\Canm\in
\CHb{\bigcup_{i\in I} \Kat A_i}$, then there exists a quotient map $\Zobr
f{\coprod\limits_{i\in I}B_i}\Canm$, where $B_i$ belongs to $\Kat A_i$ for each
$i\in I$. Put $f_i=f|_{B_i}$ and let $A_i$ be the space on the set
$\Obr{f_i}{B_i}$ endowed with the quotient topology with respect to $f_i$ for
each $i\in I$.

The topology of every space $A_i$ is finer than the topology of the
corresponding \ssp of $\Canm$ and it follows that $A_i$ is either discrete or
\filsp. Clearly, a set $U\subseteq \Canm$ is open in $\Canm$ \iaoi $U\cap A_i$
is open in $A_i$ for each $i\in I$ and $A_i\in\Kat A_i$. Obviously, there
exists $i_0\in I$ such that $\alpha$ is an accumulation point of $A_{i_0}$
(otherwise $\alpha$ would be isolated in $\Canm$).

We show that $\Canm\in \CHb{A_{i_0}}$. Let $M$ be a non-closed subset of
$\Canm$. It suffices to find a continuous map $\Zobr f{A_{i_0}}\Canm$ such that
$\Obr f{A_{i_0}\setminus\{\alpha\}} \subseteq M$ and $f(\alpha)=\alpha$.
According to Corollary \ref{CORSA} the \ssp on the set $M\cup\{\alpha\}$ is
homeomorphic to $\Canm$. Let us denote the homeomorphism from $\Canm$ to
$M\cup\{\alpha\}$ by $g$. Moreover, there is a continuous map
$\Zobr{i}{A_{i_0}}\Canm$ defined by $i(x)=x$ for each $x\in A_{i_0}$. The
desired continuous map is $f=g\circ i$. 
\end{proof}

If $X$ and $Y$ are \filsps, then a continuous map $\Zobr fXY$ is called a
\emph{prime map} if it maps only the accumulation point of $X$ to the
accumulation point of $Y$.

\begin{LM}\label{LMFG3}
Let $\alpha\geq\omega_1$ be \arcn and $0<m<n<\omnul$. There exists a quotient
prime map $\Zobr g{\Canm}{\Camm}$.
\end{LM}

\begin{proof}
Obviously, it suffices to prove the lemma for $n=m+1$. In this case
$\Camjm=P(\sum \Camm)$ is a \tsp on the set $\{\alpha\}\cup \alpha^{m+1}$ and
$\Camm$ is a \tsp on the set $\{\alpha\}\cup \alpha^m$. We define a map $\Zobr
g\Camjm\Camm$ by $g(\alpha)=\alpha$ and $g((\eta,x))=x$ for all $(\eta,x)\in
\Camjm\setminus\{\alpha\}$. It is easy to check that the map $g$ is continuous
and quotient. 
\end{proof}

\begin{COR} \label{CORSUBSET}
If $\alpha\geq\omega_1$ is \arcn and $0<m<n<\omnul$, then $\Camm\in
\CHb{\Canm}$.
\end{COR}

\begin{PROP}\label{PROPCAI0}
If $\alpha$ is \arcn and $\SCHb{\Calp}\subseteq\CHb{\Kat A\cup \Kat B}$, then
$\SCHb{\Calp}\subseteq\CHb{\Kat A}$ or $\SCHb{\Calp}\subseteq\CHb{\Kat B}$.
\end{PROP}

\begin{proof}
Since the case $\alpha=\omnul$ follows immediately from Proposition
\ref{PROPCOMI0} we can assume that $\alpha\geq\omega_1$.

By Proposition \ref{PRAI0} for each $n$, $0<n<\omnul$, the space $\Canm$
belongs either to $\Kat A$ or to $\Kat B$. By Lemma \ref{LMFG3} we have a
quotient map $\Zobr f{\Canm}{\Camm}$ for each $n>m$. Hence, one of these two
coreflective categories contains all spaces $\Canm$ and, consequently, it
contains $\SCHb{\Calp}$. 
\end{proof}

Now we can state the main result of this section.

\begin{THM}\label{THMJOIN}
If $\Kat A$, $\Kat B$ are coreflective subcategories of the category $\Top$ and
$\HCK A=\HCK B=\FG$, then $\HCKb{\CHb{\Kat A\cup\Kat B}}=\FG$.
\end{THM}

\begin{proof}
Suppose the contrary. Then according to \cite[Theorem 4.8]{CIN} there exists \arcn $\alpha$
with $\SCHb{\Calp}\subseteq\CHb{\Kat A\cup\Kat B}$. Proposition \ref{PROPCAI0} implies that
$\SCHb{\Calp}\subseteq\Kat A$ or $\SCHb{\Calp}\subseteq\Kat B$, contradicting the assumption
that the hereditary coreflective kernel of both these categories is $\FG$. 
\end{proof}

Let $\mathcal C$ be the conglomerate of all coreflective subcategories of \Top.
It is well known that $\mathcal C$ partially ordered by $\subseteq$ is a
complete lattice. Denote by $\mathcal K$ the conglomerate of all coreflective
subcategories $\Kat A$ of $\Top$ with $\HCK A=\FG$. The above theorem says that
$\mathcal K$ is closed under the formation of non-empty finite joins in
$\mathcal C$. We next show that $\mathcal K$ fails to be closed under the
formation of infinite countable joins in $\mathcal C$. In the concrete we prove
that if $\alpha\geq\omega_1$ is a regular cardinal, then all categories
$\CHb{\Canm}$ belong to $\mathcal K$. According to Proposition
\ref{PROPAMNAOMCALP} the category $\SCHb{\Calp}$ is the join of this family in
$\mathcal C$ and, evidently, $\SCHb{\Calp}\notin\mathcal K$. The proof is
divided into three auxiliary lemmas.

\begin{LM}\label{LMFG4pom}
Let $\alpha\geq\omega_1$ be \arcn and $2\leq n<\omnul$. If there exists a prime
map $\Zobr f\Canm\Canjm$, then there exists a prime map $\Zobr {f'}\Canm\Canjm$
such that $\Obr{f'}{\{\xi\}\times\alpha^{n-1}} \cap (\bigcup_{\eta<\xi}
\{\eta\}\times \alpha^{n})=\emptyset$ for each $\xi<\alpha$.
\end{LM}

\begin{proof}
Let $\Zobr f\Canm\Canjm$ be a prime map. Denote by $B_\xi$ the \ssp of $\sum \Canom$ on the
set $\{\xi\} \cup (\{\xi\} \times \alpha^{n-1})$ where $\xi<\alpha$. The \ssp $B_\xi$ is
homeomorphic to $\Canom$.

For each $\xi<\alpha$ the set $\Invobr f{\{\alpha\} \cup (\bigcup_{\eta\geq\xi}
\{\eta\} \times \alpha^n)}$ is open in $\Canm$, therefore there exists an
ordinal $\gamma<\alpha$ such that for each $\gamma'>\gamma$ the set
$\{\gamma'\} \cup (\Invobr f{\bigcup_{\eta\geq\xi} \{\eta\} \times \alpha^n}
\cap B_{\gamma'})$ is open in $B_{\gamma'}$. Hence, we can define an increasing
sequence $(\gamma_\xi)_{\xi<\alpha}$ such that $C_\xi:=\{\gamma_\xi\} \cup
(\Invobr f{\bigcup_{\eta\geq\xi} \{\eta\} \times \alpha^n} \cap
B_{\gamma_\xi})$ is open in $B_{\gamma_\xi}$. Clearly, $\Obr
f{C_\xi\setminus\{\gamma_\xi\}} \subseteq \bigcup_{\eta\geq\xi} \{\eta\} \times
\alpha^n$.

According to Corollary \ref{CORSA} the \ssp of $B_{\gamma_\xi}$ on the set
$C_\xi$ is homeomorphic to $\Canom$. Hence, for each $\xi<\alpha$ we can define
an embedding $\Vloz{h_\xi}\Canom{\sum\Canom}$ such that
$\Obr{h_\xi}\Canom=C_\xi$. It is easy to see that the map $\Zobr
h{\sum\Canom}{\sum\Canom}$ given by $h(\xi)=\gamma_\xi$ for each $\xi<\alpha$,
$h(\alpha)=\alpha$ and $h(\xi,x)=h_\xi(x)$ for each $\xi<\alpha$ and
$x\in\alpha^{n-1}$ is also an embedding. Put $A_\xi=\{\xi\}\times\alpha^{n-1}$
($A_\xi\subseteq B_\xi$). Then $\Obr {h}{A_\xi}\subseteq \Obr {h_\xi}{\Canom}
=C_\xi$ and $\Obr f{\Obr h{A_\xi}} \subseteq \Obr
f{C_\xi\setminus\{\gamma_\xi\}}\subseteq \bigcup_{\eta\geq\xi}
\{\eta\}\times\alpha^n$. Consequently, $\Obr {f\circ h}{A_\xi} \cap
(\bigcup_{\eta<\xi} \{\eta\} \times \alpha^n) = \emps$ and the prime map
$\Zobr{f'=f\circ(h|_{\Canm})}\Canm\Canjm$ is a prime map satisfying the
required condition. 
\end{proof}

\begin{LM}\label{LMFG4}
Let $\alpha\geq\omega_1$ be \arcn and $0<n<\omnul$. Then there exists no prime
map from $\Canm$ to $\Canjm$.
\end{LM}

\begin{proof}
First let $n=1$. For each $\gamma<\alpha$ the set $\{\gamma\}\times\alpha$ is
closed in $\Calp^-_2$. Consequently, $\Invobr f{\{\gamma\}\times\alpha}$ is
closed in $\Calp$, hence it contains less than $\alpha$ points and there exists
a set $U_\gamma\subseteq\alpha$ with $\card(\alpha\setminus U_\gamma)<\alpha$
such that $(\{\gamma\}\times U_\gamma) \cap \Obr f{\Calp}=\emptyset$. Thus,
$W=\{\alpha\} \cup \left( \bigcup_{\gamma<\alpha} \{\gamma\}\times U_\gamma
\right)$ is an open neighborhood of $\alpha$ in $\Calp^-_2$ such that $\Invobr
fW=\{\alpha\}$ and this contradicts the continuity of $f$.

Let $n>1$ and the lemma hold for $n-1$. Suppose that there exists a prime map
$\Zobr f\Canm\Canjm$. By Lemma \ref{LMFG4pom} we can assume \wlogg that $\Obr
f{\{\xi\}\times\alpha^{n-1}}\cap (\bigcup_{\eta<\xi} \{\eta\} \times
\alpha^{n})=\emps$ for each $\xi<\alpha$.

Recall the definition of the quotient prime map $\Zobr g{\Calp^-_n}{\Calp^-_{n-1}}$ from
Lemma \ref{LMFG3}. The map $g$ is defined by $g(\alpha)=\alpha$ and $g(\eta,x)=x$ for
$\eta<\alpha$, $x\in\alpha^{n-1}$.

Put $A_\xi=\{\xi\} \times \alpha^{n-1}$. Let us denote the \ssp of $\sum
\Canom$ on the set $\{\xi\}\cup A_\xi= \{\xi\} \cup (\{\xi\} \times
\alpha^{n-1})$ by $B_\xi$ for each $\xi<\alpha$. Clearly, $B_\xi$ is
homeomorphic to $\Canom$. We define a map $\Zobr{f_\xi}{B_\xi}\Canjm$ by
$f_\xi(\xi)=\alpha$ and $f_\xi(\xi,x)=f(\xi,x)$ for each $x\in\alpha^{n-1}$.

The map $\Zobr{g\circ f_\xi}{B_\xi}\Canm$ cannot be continuous, otherwise we
get a prime map from a space homeomorphic to $\Canom$ to the space $\Canm$.
Therefore there exists an open subset of $\Canm$ such that inverse image of
this set is not open in $B_\xi$. This set can be written in the form $U_\xi
\cup \{\alpha\}$, where $\alpha\notin U_\xi$, and we get that the set
$$\Invobr{f_\xi}{\Invobr g{U_\xi \cup \{\alpha\}}} = \Invobr{f_\xi}{\{\alpha\}
\cup (\bigcup_{\eta<\alpha} (\{\eta\}\times U_\xi))} = \{\xi\} \cup (B_\xi \cap
\Invobr f{\bigcup_{\eta<\alpha} \{\eta\}\times U_\xi})$$ is not open in
$B_\xi$.

Put $V_\xi=\bigcap_{\eta\leq\xi} U_\eta$ for $\xi<\alpha$. The family $V_\xi$
is non-increasing and it has the same properties as the family $U_\xi$. Each
$V_\xi\cup\{\alpha\}$ is open in $\Canm$, because $\Canm$ belongs to $\Topal$
($\SCHb{\Calp}\subseteq\Topal$). The set $\{\xi\} \cup (B_\xi \cap \Invobr
f{\bigcup_{\eta<\alpha} \{\eta\}\times V_\xi})$ is not open in $B_\xi$ since
$B_\xi$ is a \filsp with an accumulation point $\xi$ (and $\{\xi\} \cup (B_\xi
\cap \Invobr f{\bigcup_{\eta<\alpha} \{\eta\}\times V_\xi}) \subseteq \{\xi\}
\cup (B_\xi \cap \Invobr f{\bigcup_{\eta<\alpha} \{\eta\}\times U_\xi})$).

Finally let us put $W=\bigcup_{\xi<\alpha} \{\xi\} \times V_\xi$. The set
$W\cup\{\alpha\}$ is open in $\Canjm$. We claim that $\Invobr f{\{\alpha\}\cup
W}$ is not open in $\Canm$. It suffices to show that $\{\xi\}\cup (\Invobr
f{\{\alpha\}\cup W} \cap B_\xi)$ is not open in $B_\xi$ for each $\xi <
\alpha$.

Clearly, $B_\xi=A_\xi\cup\{\xi\}$ and we get $\{\xi\} \cup (\Invobr
f{\{\alpha\} \cup W} \cap B_\xi)= \{\xi\} \cup (\Invobr f{\bigcup_{\eta<\alpha}
\{\eta\} \times V_\eta} \cap A_\xi)$. We have $\Obr f{A_\xi} \cap
(\bigcup_{\eta<\xi} \{\eta\} \times \alpha^{n-1})=\emps$, hence $\Invobr
f{\bigcup_{\eta<\alpha} \{\eta\} \times V_\eta} \cap A_\xi = \Invobr
f{\bigcup_{\eta\geq\xi} \{\eta\} \times V_\eta} \cap A_\xi$ and we obtain
\begin{gather*}
\{\xi\} \cup (\Invobr f{\{\alpha\} \cup W} \cap B_\xi) =
\{\xi\} \cup (\Invobr f{\bigcup_{\eta\geq\xi} \{\eta\}\times V_\eta} \cap B_\xi) \subseteq \\
\subseteq \{\xi\} \cup (\Invobr f{\bigcup_{\eta\geq\xi} \{\eta\}\times V_\xi} \cap B_\xi)
\subseteq \{\xi\} \cup (\Invobr f{\bigcup_{\eta<\alpha} \{\eta\}\times V_\xi}
\cap B_\xi).
\end{gather*}
The latter set is not open in $B_\xi$ therefore $\{\xi\} \cup (\Invobr
f{\{\alpha\}\cup W} \cap B_\xi)$ is not open in $B_\xi$ as well. 
\end{proof}

\begin{LM}\label{LMHCKCALPN}
Let $\alpha\geq\omega_1$ be \arcn and $0<n<\omnul$. Then
$\HCKb{\CHb{\Canm}}=\FG$.
\end{LM}

\begin{proof}
Recall (see \cite{HER}) that if $\gamma>\delta$, then $\Top(\gamma) \cap
\Gen(\delta)=\FG$. For $\beta<\alpha$ we have $\SCHb{\Cbet}\subseteq
\Gen(\beta)$ and $\Canm \in \Topal$, hence $\SCHb{\Cbet}\nsubseteq
\CHb{\Canm}$. Similarly if $\beta>\alpha$, then $\SCHb{\Cbet}\subseteq
\Top(\beta)$ and $\Canm \in \Genal$. Thus, $\SCHb{\Cbet}\nsubseteq
\CHb{\Canm}$.

By Lemma \ref{LMFG4} and Lemma \ref{LMNOVA2} $\Canjm \notin \CHb{\Canm}$ (every
\filssp of $\Canm$ is homeomorphic to $\Canm$) and $\Canjm\in\SCHb{\Calp}$,
therefore $\SCHb{\Calp} \nsubseteq \CHb{\Canjm}$ as well. 
\end{proof}

Denote by $\mathcal L$ the collection of all coreflective subcategories $\Kat
A$ of $\Top$ such that $\SSp A=\Top$ and $\HCK A=\FG$. In the paper \cite{SLEZ}
it is shown that $\mathcal L$ has the smallest element $\Kat
A_0=\CHb{\{S^\alpha; \alpha$ is a cardinal$\}}$, where $S$ is the \Sierp
doubleton, and $\mathcal L$ is closed under the formation of arbitrary
non-empty intersections. This together with Theorem \ref{THMJOIN} yields:

\begin{THM}
The collection $\mathcal L$ is closed under the formation of non-empty
intersections, non-empty finite joins in $\mathcal C$ and has the smallest
element.
\end{THM}

\begin{PROP}
There is no maximal coreflective subcategory $\Kat A$ of $\Top$
such that $\HCK A=\FG$. Consequently, the collection $\mathcal L$
has no maximal element.
\end{PROP}

\begin{proof}
Suppose that $\Kat A$ is maximal coreflective subcategory of $\Top$ with the
property $\HCK A=\FG$. Let $\alpha\geq\omega_1$ be a regular cardinal.
According to Lemma \ref{LMHCKCALPN} and Theorem \ref{THMJOIN} $\HCKb{\CHb{\Kat
A\cup \{\Canm\}}}=\FG$ for each $n$, $0<n<\omnul$. Thus, we get $\Canm\in\Kat
A$ for each $n$ and by Proposition \ref{PROPAMNAOMCALP}
$\SCHb\Calp\subseteq\Kat A$, a contradiction.

The proof that $\mathcal L$ has no maximal elements is analogous.
\end{proof}

The family $\CHb{\Kat A_0\cup \{\Canm\}}$, $0<n<\omnul$, where
$\alpha\geq\omega_1$ is \arcn, is an example of a countable family of elements
of $\mathcal L$ such that its join does not belong to $\mathcal L$.

\noindent\textbf{Acknowledgement.}
The author wishes to thank J.~\v{C}in\v{c}ura for many helpful comments made on
earlier drafts of this paper.

\end{document}